\documentclass[a4paper,11pt]{article}
\usepackage{amsmath}
\usepackage{amssymb}
\usepackage{enumerate}
\usepackage{theorem}
\usepackage{array}
\usepackage{bm}
\usepackage{mathrsfs}
\newtheorem{theorem}{Theorem}[section]
\newtheorem{lemma}[theorem]{Lemma}
\newtheorem{corollary}[theorem]{Corollary}
\newtheorem{proposition}[theorem]{Proposition}
\theorembodyfont{\rmfamily}
\newtheorem{definition}[theorem]{Definition}
\newtheorem{notation}[theorem]{Notation}
\newtheorem{remark}[theorem]{Remark}

\newcommand{\proof}{\noindent \mbox{\em Proof.\hspace*{2mm}}}
\newcommand{\qed}{\hfill \mbox{$  \Box $}}

\DeclareMathOperator{\image}{Im}

\DeclareMathOperator{\Fix}{Fix}

\title{On maximal plane curves of degree $3$ over $\mathbb{F}_4$,
and Sziklai's example of degree $q-1$ over $\mathbb{F}_q$}
\author{
Masaaki Homma
\\
Department of Mathematics and Physics\\
Kanagawa University\\
Hiratsuka 259-1293, Japan\\
homma@kanagawa-u.ac.jp
}

\date{}

\begin{document}
\maketitle
\begin{abstract}
The classification of maximal plane curves of degree $3$
over $\mathbb{F}_4$ will be given, which complements 
Hirschfeld-Storme-Thas-Voloch's theorem on a characterization of Hermitian curves in $\mathbb{P}^2$.
This complementary part should be understood as the classification of Sziklai's example of maximal plane curves of degree $q-1$ over $\mathbb{F}_q$.
Although two maximal plane curves of degree $3$ over $\mathbb{F}_4$
up  to projective equivalence over $\mathbb{F}_4$ appear,
they are birationally equivalent over $\mathbb{F}_4$ each other.
\end{abstract}

\section{Introduction}
This paper is concerned with upper bounds for the number of $\mathbb{F}_q$-points of plane curves defined over $\mathbb{F}_q$.
Let $C$ be a plane curve defined by a homogeneous equation $f \in \mathbb{F}_q[x_0, x_1, x_2]$.
The set of $\mathbb{F}_q$-points $C(\mathbb{F}_q)$ of $C$ is
$\{ (a_0, a_1, a_2) \in \mathbb{P}^2 \mid a_0, a_1, a_2 \in \mathbb{F}_q \text{ and } f(a_0, a_1, a_2) =0\}$.
The cardinality of $C(\mathbb{F}_q)$ is denoted by $N_q(C)$, and the degree of $C$ by $\deg C$, or simply by $d$.
We are interesting in upper bounds for $N_q(C)$ with respect to $\deg C$.

Aubry-Perret's generalization \cite{aub-per1996} of the Hasse-Weil bound implies that for absolutely irreducible plane curve of degree $d$ over $\mathbb{F}_q$,
\begin{equation}\label{APHW}
N_q(C) \leq q+1 +(d-1)(d-2)\sqrt{q}.
\end{equation}
On the other hand, the Sziklai bound established by a series of papers of Kim and the author \cite{hom-kim2009, hom-kim2010a, hom-kim2010b} gives one under a more mild condition, that is, for $C$ without $\mathbb{F}_q$-linear components,
\begin{equation}\label{Sziklai}
N_q(C) \leq (d-1)q +1
\end{equation}
except for the curve over $\mathbb{F}_4$ defined by
\[
(x_0+x_1+x_2)^4+(x_0x_1+x_1x_2+x_2x_0)^2+x_0x_1x_2(x_0+x_1+x_2)=0.
\]
When $d < \sqrt{q}+1$, the Aubry-Perret generalization of Hasse-Weil bound is better than the Sziklai bound, however when $d > \sqrt{q} +1$, the latter is better than the former, and these two bounds meet at $d=\sqrt{q}+1$, that is,
both (\ref{APHW}) and (\ref{Sziklai}) imply
\begin{equation}\label{qplusone}
N_q(C) \leq \sqrt{q}^3 + 1 \text{ if } \deg C = \sqrt{q}+1,
\end{equation}
where $q$ is an even power of a prime number. 
From now on, when a statement contains $\sqrt{q}$, we tacitly understand $q$ to be an even power of a prime number. 

\medskip

Three decades ago, 
Hirschfeld, Storme, Thas and Voloch \cite{hir-sto-tha-vol1991}
gave a characterization of Hermitian curves of degree $\sqrt{q}+1$
over $\mathbb{F}_q$, which is a maximal curve in the sense of the bound (\ref{qplusone}).
Here we understand a Hermitian curve as a plane curve defined by an equation
\[
(x_0^{\sqrt{q}}, x_1^{\sqrt{q}}, x_2^{\sqrt{q}}) A
\begin{pmatrix}
x_0\\  x_1\\  x_2
\end{pmatrix}
=0
\]
for a certain matrix $A \in GL(3, \mathbb{F}_q)$
satisfying ${}^t\! A =A^{(\sqrt{q})}$,
where
${}^t\! A$ denotes the transposed matrix of $A$
and 
$A^{(\sqrt{q})}$ the matrix taking entry-wise the $\sqrt{q}$-th power of $A$.
Note that any two Hermitian curves are projectively equivalent each other
over $\mathbb{F}_q$ \cite[\S 7.3]{hir}.
\begin{theorem}[Hirschfeld-Storme-Thas-Voloch]
In $\mathbb{P}^2$ over $\mathbb{F}_q$ with $q \neq 4$, a curve over $\mathbb{F}_q$ of degree $\sqrt{q}+1$, without $\mathbb{F}_q$-linear components,
which contains $\sqrt{q}^3+1$ $\mathbb{F}_q$-points, is a Hermitian curve.
\end{theorem}
For $q=4$, they gave an example of a nonsingular plane curve over $\mathbb{F}_4$ which had $9 \,(= 2^3 +1)$ $\mathbb{F}_4$-points,
but was not a Hermitian curve.
Actually the plane curve defined by
\begin{equation}\label{counter_example}
x_0^3 + \omega x_1^3 + \omega^2 x_2^3 =0
\end{equation}
is such an example, where
$\mathbb{F}_4 =\{0, 1, \omega, \omega^2  \}.$

Our primary concern is to complete the determination of plane curves over $\mathbb{F}_q$ of degree $\sqrt{q}+1$ with $\sqrt{q}^3+1$ $\mathbb{F}_q$-points.

\begin{theorem}\label{maintheorem}
Let $C$ be a plane curve over $\mathbb{F}_q$ without $\mathbb{F}_q$-linear components. If $\deg C =\sqrt{q}+1$ and $N_q(C)= \sqrt{q}^3 +1$, then $C$ is either
\begin{enumerate}[{\rm (i)}]
\item a Hermitian curve, or
\item a nonsingular curve of degree $3$ which is projectively equivalent to the curve {\rm (\ref{counter_example})} over $\mathbb{F}_4$.
\end{enumerate}
\end{theorem}

The second case (ii) in the above theorem should be understood the case of $q=4$ among Sziklai curves \cite{szi2008} of degree $q-1$ that achieve the Sziklai bound  (\ref{Sziklai}).
Here a Sziklai curve means one over $\mathbb{F}_q$, of degree $q-1$
defined by the following type of equation:
\begin{equation}\label{Sziklai_example}
\alpha x_0^{q-1} + \beta x_1^{q-1} + \gamma x_2^{q-1}=0 
\ \text{with}\  \alpha\beta\gamma \neq 0 \ \text{and}\  \alpha + \beta + \gamma=0.
\end{equation}

The curve (\ref{Sziklai_example}) will be denoted by $C_{(\alpha, \beta, \gamma)}$.
Since $x^{q-1}=1$ for any $x\in \mathbb{F}_q^\ast$ and $\alpha + \beta + \gamma=1$,
\begin{equation}\label{only_supset}
C_{(\alpha, \beta, \gamma)}(\mathbb{F}_q) \supset
\mathbb{P}^2(\mathbb{F}_q) \setminus (\cup_{i=0}^{2} \{ x_i =0\}).
\end{equation}
Here $\{ x_i =0\}$ denotes the line defined by $x_i=0$.
Furthermore, since $\deg C_{(\alpha, \beta, \gamma)} =q-1$,
\[
N_q(C_{(\alpha, \beta, \gamma)}) \leq (q-2)q +1 = (q-1)^2
\]
by the Szikali bound. Therefore equality must hold in (\ref{only_supset}),
that is,
\begin{equation}\label{rational_Sziklai}
C_{(\alpha, \beta, \gamma)}(\mathbb{F}_q)=
\mathbb{P}^2(\mathbb{F}_q) \setminus (\{ x_0 =0\}\cup\{ x_1 =0\}\cup\{ x_2 =0\}).
\end{equation}
Note that $C_{(\alpha, \beta, \gamma)}$ makes sense under the condition $q>2$.

\begin{theorem}\label{classification}
The number $\nu_q$ of projective equivalence classes over $\mathbb{F}_q$ in the family of curves
\[
 \{C_{(\alpha, \beta, \gamma)} \mid \alpha, \beta , \gamma \in \mathbb{F}_q^{\ast}, \ \alpha + \beta + \gamma = 0 \}
\]
is as follows:
\begin{enumerate}[{\rm (I)}]
\item Suppose that the characteristic of $\mathbb{F}_q$ is neither $2$ nor $3$.
 \begin{enumerate}[{\rm (\mbox{I}-i)}]
  \item If $q \equiv 2 \bmod 3$, then $\nu_q = \frac{q+1}{6}$.
  \item If $q \equiv 1 \bmod 3$, then $\nu_q = \frac{q+5}{6}$.
 \end{enumerate}
\item Suppose that $q$ is a power of $3$. Then  $\nu_q = \frac{q+3}{6}$.
\item Suppose that  $q$ is a power of $2$.
 \begin{enumerate}[{\rm (\mbox{III}-i)}]
  \item If $q =2^{2s+1}$, that is, $q\equiv 2 \bmod 3$, then $\nu_q = \frac{q-2}{6}$.
  \item If $q =2^{2s}$, that is, $q\equiv 1 \bmod 3$, then $\nu_q = \frac{q+2}{6}$.
 \end{enumerate}
\end{enumerate}
\end{theorem}
In this theorem, we don't assume $q>2$ explicitly, however the assertion (III-i) says the family of curves in question is empty if $q=2$.

\medskip

The construction of this article is as follows:

In Section~2, we will give the proof of Theorem~\ref{classification}
together with the characterization of Sziklai curves of degree $q-1$.

In Section~3, we will give the proof of Theorem~\ref{maintheorem};
actually we will handle the case $q=4$.

In Section~4, we will make explicitly an $\mathbb{F}_4$-isomorphism between the function field of the Hermitian curve over $\mathbb{F}_4$ defined by $x_0^3 + x_1^3 + x_2^3=0$ and that of the curve (\ref{counter_example}).

\section{Sziklai's example of maximal curves of degree $q-1$}
The purpose of this section is to prove Theorem~\ref{classification}.
Let 
$
\mathscr{S}_q =
\{C_{(\alpha, \beta, \gamma)} \mid \alpha, \beta , \gamma \in \mathbb{F}_q^{\ast}, \ \alpha + \beta + \gamma = 0 \}.
$
The first step of the proof is to give a characterization of the member of $\mathscr{S}_q$.

\begin{proposition}\label{prop_characterization}
Let $C$ be a possibly reducible plane curve over $\mathbb{F}_q$
of degree $q-1$.
Then $C \in \mathscr{S}_q$
if and only if
\begin{equation}\label{characterization}
C(\mathbb{F}_q) = \mathbb{P}^2(\mathbb{F}_q) \setminus \left( \bigcup_{i=0}^2 \{x_i = 0\}\right).
\end{equation}
\end{proposition}
The ``only if" part has already observed in Introduction.
Now we prove the ``if" part.
\begin{lemma}\label{lemma_affine}
In $\mathbb{A}^2$ with coordinates $x, y$ over $\mathbb{F}_q$,
the ideal $I$ in $\mathbb{F}_q[x,y]$ of the set
$
\{ (a,b) \in \mathbb{F}_q^2 \mid ab \neq 0 \}
$
is $(x^{q-1}-1, y^{q-1}-1).$

Furthermore, if $f(x,y) \in I$ is of degree at most $q-1$,
then $f(x,y) = \alpha (x^{q-1} -1) + \beta (y^{q-1}-1)$ for some $\alpha , \beta \in \mathbb{F}_q$.
\end{lemma}
\proof
Let $J$ denote the ideal $(x^{q-1}-1, y^{q-1}-1)$ of $\mathbb{F}_q[x,y]$.
Obviously $J \subseteq I$.
For $f(x,y) \in I$,
there are polynomials $g_i(x) \in  \mathbb{F}_q[x]$ ($0\leq i \leq q-2$)
of degree $\leq q-2$ so that
\[
f(x,y) \equiv \sum_{i=0}^{q-2} g_i(x)y^i  \bmod J.
\]
For each $a \in \mathbb{F}_q^{\ast}$,
the equation $\sum_{i=0}^{q-2} g_i(a)y^i=0$ has to have
$q-1 \, (= |\mathbb{F}_q^{\ast}|)$ solutions
because  $\sum_{i=0}^{q-2} g_i(x)y^i \in I$.
Hence $g_i(a)=0$ for any $i$.
Since $\deg g_i \leq q-2$,
$g_i$ must be the zero polynomial.
Hence $f(x,y) \equiv  0 \bmod J$.
This completes the proof of the first part.

For the second part, let $\alpha$ and $\beta$ be the coefficients of $x^{q-1}$ and $y^{q-1}$ in $f(x,y)$ respectively.
Then
\begin{equation}\label{qminusone}
f(x,y)-\alpha(x^{q-1}-1) -\beta(y^{q-1}-1)
= \sum_{i=1}^{q-2} u_{q-1-i}(x)y^i + v_{q-2}(x),
\end{equation}
where $\deg u_{q-1-i}(x) \leq q-1-i \ (\leq q-2)$ and
$\deg v_{q-2}(x) \leq q-2$.
So the same argument as above works well,
and we know the right side of (\ref{qminusone})
is the zero polynomial.
\qed 

\medskip

\noindent \mbox{\em Proof of Proposition}~\ref{prop_characterization}.\hspace*{2mm}
Choose a homogeneous equation $f(x_0, x_1, x_2)=0$ of degree $q-1$ over $\mathbb{F}_q$ for a given curve $C$ with the property (\ref{characterization}).
From Lemma~\ref{lemma_affine},
there are elements $\alpha, \beta \in \mathbb{F}_q$
such that
$f(\frac{x_0}{x_2}, \frac{x_1}{x_2}, 1) = \alpha ((\frac{x_0}{x_2})^{q-1}-1)
 + \beta ((\frac{x_1}{x_2})^{q-1}-1).$
Therefore
$
f(x_0, x_1, x_2) = x_2^{q-1}f(\frac{x_0}{x_2}, \frac{x_1}{x_2}, 1)=
\alpha (x_0^{q-1}-x_2^{q-1}) + \beta (x_1^{q-1}-x_2^{q-1}).
$ Since $C(\mathbb{F}_q) \cap \{x_2 = 0\}$ is empty,
$f(a,b,0 ) \neq 0$ for any $(a,b) \in \mathbb{F}_q^2 \setminus \{(0,0)\}.$
In particular, $\alpha =f(1,0,0) \neq 0$,
$\beta = f(0,1,0) \neq 0$ and
$\alpha + \beta = f(1,1,0) \neq 0$.
Hence $C \in \mathscr{S}_q$.
\qed

\medskip

Now we want to classify $\mathscr{S}_q$ up to projective equivalence over $\mathbb{F}_q$.

\begin{definition}
Let $C$ be a possibly reducible curve in $\mathbb{P}^2$ over $\mathbb{F}_q$,
and $\delta$ a nonnegative integer.
An $\mathbb{F}_q$-line $l$ is said to be a $\delta$-line with respect to $C$
if $|l \cap C(\mathbb{F}_q)|=\delta$.
\end{definition}

\begin{lemma}\label{deltaline}
Let $C \in \mathscr{S}_q$, and $\delta$ a nonnegative integer such that
a $\delta$-line with respect to $C$ actually exists.
Then $\delta$ is either $0$ or $q-2$ or $q-1$, and the number of  $\delta$-lines are as in Table~{\rm \ref{delta_line}}.
\begin{table}[hbtp]
\centering
\begin{tabular}{c|c}
$\delta$ & the number of $\delta$-lines \\
\hline
$0$ & $3$ \\
$q-2$& $(q-1)^2$ \\
$q-1$ & $3(q-1)$
\end{tabular}
\caption{$\delta$-lines w.r.t. $C \in \mathscr{S}_q$}
\label{delta_line}
\end{table}
\end{lemma}
\proof
Note that $q>2$ because $\mathscr{S}_q$ is not empty.
Since $\mathbb{P}^2(\mathbb{F}_q) = C(\mathbb{F}_q) \sqcup (\cup_{i=0}^2 \{x_i =0\})$ (where the symbol $\sqcup$ indicates disjoint union) and $q>2$, the possible values of $\delta$ are $0$, $q-2$ and $q-1$.
Obviously the number of $0$-lines is $3$. A $(q-1)$-line is not a $0$-line, and passes through one of intersection points of two $0$-lines.
Other lines are $(q-2)$-lines.
\qed

\medskip

We need an elementary fact on the finite group action, so called ``Burnside's lemma" \cite[Corollary 7.2.9]{ste2012}.
\begin{lemma}\label{Burnside}
Let $G$ be a finite group which acts on a finite set $X$. For $g \in G$,
$\Fix g$ denotes the set of fixed points of $g$ on $X$.
Then the number $\nu$ of orbits of $G$ on $X$ is given by
\[
\nu = \frac{1}{|G|} \sum_{g \in G}\, |\Fix g|.
\]
\end{lemma}

\medskip

\noindent \mbox{\em Proof of Theorem}~\ref{classification}.\hspace*{2mm}
The first claim is that if two members $C_{(\alpha, \beta, \gamma)},
C_{(\alpha', \beta', \gamma')} \in \mathscr{S}_q$ are projectively equivalent over $\mathbb{F}_q$, then the point $(\alpha', \beta', \gamma') \in \mathbb{P}^2(\mathbb{F}_q)$ is a permutation of the point $(\alpha, \beta, \gamma)\in \mathbb{P}^2(\mathbb{F}_q)$, that is, there is a nonzero element $\lambda \in \mathbb{F}_q^{\ast}$ such that the triple
$(\lambda\alpha', \lambda\beta',\lambda\gamma')$
is a permutation of the triple $(\alpha, \beta, \gamma)$.

Actually, let $\Sigma$ be a projective transformation
so that $\Sigma C_{(\alpha, \beta, \gamma)}= C_{(\alpha', \beta', \gamma')}$.
Note that $\Sigma$ induces an automorphism of the homogeneous coordinate ring
$\mathbb{F}_q[x_0, x_1, x_2]$, which is denoted by $\Sigma^{\ast}$.
The set of $0$-lines with respect to each of curves in $\mathscr{S}_q$ is
$\{ \{x_0=0\},  \{x_1=0\},  \{x_2=0\}\}$
by Lemma~\ref{deltaline}.
Hence $\Sigma$ induces a permutation of those three lines.
Hence $\Sigma^{\ast}(x_i) = u_ix_{\sigma(i)}$ for some $u_i \in \mathbb{F}_q^{\ast}$, and $(\sigma(0), \sigma(1), \sigma(2))$ is a permitation of $(0,1,2)$.
Hence
\[
\Sigma^{\ast}(\alpha x_0^{q-1} + \beta x_1^{q-1} + \gamma x_2^{q-1})
= \alpha x_{\sigma(0)}^{q-1} + \beta x_{\sigma(1)}^{q-1} + \gamma x_{\sigma(2)}^{q-1}
\]
because $u_i^{q-1}=1$.

So we need to classfy $\mathscr{S}_q/ \mathbb{F}_q^{\ast}$ by the action of $S_3$ as permutations
 on coefficients. 

Observe the map
\[
\rho : \mathscr{S}_q/ \mathbb{F}_q^{\ast} \ni C_{(\alpha, \beta , \gamma)} \to
 (\alpha : \beta) \in \mathbb{P}^1(\mathbb{F}_q),
\]
which is well-defined and 
\[
\image \rho = \mathbb{P}^1(\mathbb{F}_q) \setminus \{ (0,1), (1,0), (1, -1)\}.
\]
Obviously, $\rho$ gives a one to one correspondence, so $S_3$ acts on $\image \rho$ also. Table~\ref{action} shows the $S_3$-action on $\image \rho$ explicitly.\begin{table}[hbtp]
\centering
\begin{tabular}{c|c|c}
$S_3$ & $\mathscr{S}_q/ \mathbb{F}_q^{\ast}$ & $\image \rho$ \\
\hline
$(1)$ &  $(\alpha, \beta, \gamma) \mapsto (\alpha, \beta, \gamma)$ 
& $(\alpha : \beta) \mapsto (\alpha : \beta)$\\
$(1,2)$ &  $(\alpha,  \beta, \gamma) \mapsto (\beta, \alpha, \gamma)$ 
& $(\alpha : \beta) \mapsto (\beta : \alpha)$\\
$(2,3)$ &  $(\alpha,  \beta, \gamma) \mapsto ( \alpha, \gamma, \beta)$ 
& $(\alpha : \beta) \mapsto ( \alpha :-(\alpha +\beta))$\\
$(1,3)$ &  $(\alpha,  \beta, \gamma) \mapsto ( \gamma, \beta, \alpha)$ 
& $(\alpha : \beta) \mapsto (-(\alpha +\beta) : \beta)$\\
$(1, 2,3)$ &  $(\alpha,  \beta, \gamma) \mapsto (  \gamma, \alpha,\beta)$ 
& $(\alpha : \beta) \mapsto (-(\alpha +\beta) :  \alpha )$\\
$(1, 3, 2)$ &  $(\alpha,  \beta, \gamma) \mapsto ( \beta, \gamma, \alpha)$ 
& $(\alpha : \beta) \mapsto ( \beta :-(\alpha +\beta))$
\end{tabular}
\caption{$S_3$-action on $\image \rho$}
\label{action}
\end{table}

Now we compute the number of fixed points on $\image \rho$ 
by each $\sigma \in S_3$.
\begin{itemize}
\item Fixed points of the identity $(1)$ are all the $q-2$ points of $\image \rho$.
\item $(\alpha :\beta) \in \Fix(1,2) \Leftrightarrow (\alpha :\beta) = (\beta : \alpha)  \Leftrightarrow \alpha^2-\beta^2 =0$.
If the characteristic of $\mathbb{F}_q \neq 2$, then $\Fix(1,2) = \{(1:1)\}$
because $(1:-1) \not\in \image \rho$.
If $q$ is a power of $2$, then  $\Fix(1,2)$ is empty.
\item $(\alpha :\beta) \in \Fix(2,3) \Leftrightarrow (\alpha : \beta) = ( \alpha :-(\alpha +\beta)) \Leftrightarrow \alpha = -2\beta \text{ because $\alpha \neq 0$ }.$
If the characteristic of $\mathbb{F}_q \neq 2$, then $\Fix(2,3) = \{(-2:1)\}$.
If $q$ is a power of $2$, then  $\Fix(2,3)$ is empty.
\item $(\alpha :\beta) \in \Fix(1,3) \Leftrightarrow (\alpha : \beta) = (-(\alpha +\beta) : \beta) \Leftrightarrow \beta = - 2\alpha  \text{ because $\beta \neq 0$ }.$ 
If the characteristic of $\mathbb{F}_q \neq 2$, then $\Fix(1,3) = \{(1:-2)\}$.
If $q$ is a power of $2$, then  $\Fix(1,3)$ is empty.
\item $(\alpha :\beta) \in \Fix(1,2, 3) \Leftrightarrow (\alpha : \beta) = (-(\alpha +\beta) : \alpha) \Leftrightarrow \alpha^2+ \alpha\beta +\beta^2=0 
\Leftrightarrow (\alpha : \beta) = (\eta: 1) \text{ with } \eta^2 + \eta +1 =0 \text{ and } \eta \in \mathbb{F}_q. $ 
\item $(\alpha :\beta) \in \Fix(1,3, 2) \Leftrightarrow (\alpha : \beta) = (\beta :-(\alpha +\beta)) \Leftrightarrow \alpha^2+ \alpha\beta +\beta^2=0 
\Leftrightarrow (\alpha : \beta) = (\eta: 1) \text{ with } \eta^2 + \eta +1 =0 \text{ and } \eta \in \mathbb{F}_q. $ 
\end{itemize}
For the last two cases, since a cubic root of $1$ other than $1$ exists in $\mathbb{F}_q$ if and only if $q \equiv 1 \bmod 3$,
and only the cubic root of $1$ is 1 if $q$ is a power of $3$,
\[
|\Fix (1,2,3)| = |\Fix (1,3,2)| = 
\begin{cases}
2 & \text{if $q \equiv 1 \bmod 3$}\\
1 & \text{if $q$ is a power of $3$}\\
0 & \text{else}.
\end{cases}
\]
The number of fixed points can be summarized as in Table~\ref{fixedpoints}.
\begin{table}[hbtp]
\centering
\begin{tabular}{c||c|c|c|c|c|c}
Case & $|\Fix (1)|$& $|\Fix (12)|$&$|\Fix (13)|$&$|\Fix (23)|$&
$|\Fix (123)|$&$|\Fix (132)|$\\
\hline
(I-i) & $q-2$ & $1$ & $1$ & $1$ & $0$ & $0$ \\
(I-ii)&$q-2$ & $1$ & $1$ & $1$ & $2$ & $2$ \\
(II) &$q-2$ & $1$ & $1$ & $1$ & $1$ & $1$ \\
(III-i)&$q-2$ & $0$ & $0$ & $0$ & $0$ & $0$ \\
(III-ii)&$q-2$ & $0$ & $0$ & $0$ & $2$ & $2$
\end{tabular}
\caption{Number of fixed points}
\label{fixedpoints}
\end{table}

Since $ \nu_q = \frac{1}{6}\sum_{\sigma \in S_3} |\Fix \sigma |$ by Lemma~\ref{Burnside}, we are able to know $\nu_q$ explicitly.
\qed

\medskip

At the end of this section, we raise a question:
are there non-Sziklai curves over $\mathbb{F}_q$ of degree $q-1$ that attain the Sziklai bound (\ref{Sziklai})?

\section{Maximal curves of degree $3$ over $\mathbb{F}_4$}
Let $C$ be a plane curve of degree $3$ over $\mathbb{F}_4$ without
 $\mathbb{F}_4$-linear components, and $N_4(C)=9$.
Since the degree of $C$ is $3$, $C$ is absolutely irreducible.
If $C$ had a singular point, then $C$ would be an image of $\mathbb{P}^1$
with exactly one singular point, and hence $N_4(C)$ would be at most $6\, (= N_4(\mathbb{P}^1) + 1)$.
Therefore $C$ is nonsingular.

Thanks to the Hirschfeld-Storme-Thas-Voloch theorem,
only the missing case for the classification of maximal curves of degree
$\sqrt{q}+1$ is the case of $q=4$.

\begin{theorem}\label{primarytheorem}
Let $C$ be a nonsingular plane curve of degree $3$ over $\mathbb{F}_4$.
If $N_4(C)=9$, then $C$ is either
\begin{enumerate}[{\rm(i)}]
\item Hermitian, or
\item projectively equivalent to the curve
\[
x_0^3 + \omega x_1^3 + \omega^2 x_2^3 = 0,
\]
where $\mathbb{F}_4 =\{0, 1, \omega, \omega^2  \}.$
\end{enumerate}
\end{theorem}

\begin{notation}
Let $l$ be an $\mathbb{F}_4$-line in $\mathbb{P}^2$.
The symbol $l.C$ denotes the divisor
$\sum_{P\in l \cap C} i(l.C;P)P$ on $C$,
where $i(l.C;P)$ is the local intersection multiplicity of $l$ and $C$ at $P$.
Note that though $l.C$ is defined over $\mathbb{F}_4$, a point $P$ in the support of $l.C$ may not be $\mathbb{F}_4$-point.
\end{notation}

From now on, we consider a nonsingular plane curve $C$ of degree $3$ with $N_4(C)=9$,
and lines over $\mathbb{F}_4$.

\begin{lemma}\label{ordinarytangent}
Let $l$ be a $2$-line with respect to $C$, say
$ l\cap C(\mathbb{F}_4) =\{ P_1, P_2\}$.
Then
$l.C= 2P_1+P_2$ or $P_1+2P_2$.
\end{lemma}
\proof
Since $\deg C=3$, there is a closed point $Q$ of $C$
such that $l.C = P_1 + P_2 +Q$.
Applying the Frobenius map $F_4$  over $\mathbb{F}_4$
to both side of the above equality, we know 
$P_1 + P_2 +Q= P_1 + P_2 +F_4(Q)$,
which implies that the point $Q$ is also $\mathbb{F}_4$-point.
Therefore $Q$ must concide with either $P_1$ or $P_2$ because
$l$ is a 2-line.
\qed

\begin{lemma}\label{flex}
Let $l_0$ be a $1$-line with respect to $C$, say $l_0 \cap C(\mathbb{F}_4)= \{P\}$.
Then $l_0 . C = 3P_0$.
\end{lemma}
\proof 
Consider all the $\mathbb{F}_4$-lines passing through the point $P$,
say $l_0, l_1, \dots , l_4$.
Counting $N_4(C)$ by using the disjoint union
\[
C(\mathbb{F}_q)
= \{P\} \sqcup \left(
\sqcup_{i=1}^{4}( l_i\cap C(\mathbb{F}_4) \setminus \{P \})
\right),
\]
we know that  $|l_i\cap C(\mathbb{F}_4) \setminus \{P \}|$ is $2$,
that is the remaining four lines $l_1 , \dots l_4$ to be
$3$-lines with respect to $C$.
So  each of them meets with $C$ transversally because $\deg C =3$.
Therefore $l_0$ is the tangent line to $C$ at $P$. Hence
there is a closed point $Q \in C$ such
that $l_0.C =2P +Q$. Apply $F_4$ to this divisor, $Q$ should be $\mathbb{F}_4$-points. Since $l_0$ is a $1$-line, $Q=P$.
\qed

\begin{definition}
Since $C$ is nonsingular, for any closed point $P\in C$,
the tangent line to $C$ at $P$ exists, which is a unique line $l$ such that $i(l.C;P)\geq 2$. This line is denoted by $T_P(C)$. A point $P$ with
$i(T_P(C).C;P)= 3$ is called a flex or an inflection point.
It is obvious that if $P$ is an $\mathbb{F}_4$-points, then $T_P(C)$ is an $\mathbb{F}_4$-line.
\end{definition}

\begin{corollary}\label{tangentline}
Let $P \in C(\mathbb{F}_4).$
\begin{enumerate}[{\rm (i)}]
\item If $i(T_P(C).C;P) =3$, then $T_P(C)$ is a $1$-line,
 and conversely, if an $\mathbb{F}_4$-line $l$ passing through $P$ is a $1$-line,
 then $l= T_P(C)$ and $i(T_P(C).C;P) =3$.
\item If $i(T_P(C).C;P) =2$, then $T_P(C)$ is a $2$-line, and conversely,
if an $\mathbb{F}_4$-line $l$ passing through $P_1, P_2 \in C(\mathbb{F}_4)$ is a $2$-line, then $l$ coincides with either $T_{P_1}(C)$ or $T_{P_2}(C)$.
\end{enumerate}
\end{corollary}
\proof
(i) The first part is obvious because $\deg C=3$,
and the second part is a consequence of Lemma~\ref{flex}.

(ii) This is also a consequence of Lemma~\ref{flex}:
since $T_P(C)$ is not a $1$-line, it should be a $2$-line, and the second part
is just in Lemma~\ref{ordinarytangent}
\qed

\begin{notation}
For each nonnegative integer $\delta \leq 3$,
$\mathscr{L}_{\delta}$ denotes the set of $\delta$-lines
with respect to $C$, and
$\mu_{\delta}$ denotes the cardinality of the set $\mathscr{L}_{\delta}$.
\end{notation}

The next lemma is essential for the proof of Theorem~\ref{primarytheorem}.

\begin{lemma}\label{essentiallemma}
The possibilities of quadruple $(\mu_0, \mu_1, \mu_2, \mu_3)$
are either
\begin{enumerate}[{\rm (i)}]
\item $\mu_0=0,\  \mu_1=9, \ \mu_2=0,\  \mu_3=12$; or
\item $\mu_0=3,\  \mu_1=0, \ \mu_2=9,\  \mu_3=9$.
\end{enumerate}
\end{lemma}
\proof
{\em Step}~1.
Let us consider the correspondence
\[
\mathscr{I}:=\{ (l, P) \in \breve{\mathbb{P}}^2(\mathbb{F}_4) \times C(\mathbb{F}_4) \mid l \ni P\}
\]
with projections $p_1: \mathscr{I} \to \breve{\mathbb{P}}^2(\mathbb{F}_4)$
and $p_2: \mathscr{I} \to C(\mathbb{F}_4)$,
where $\breve{\mathbb{P}}^2(\mathbb{F}_4)$ is the projective space of the $\mathbb{F}_4$-lines.
Since $|p_2^{-1}(P)| = 5$ for all $P \in C(\mathbb{F}_4)$ and $|C(\mathbb{F}_4)|= 9$, we know $|\mathscr{I}|=45$.

From Corollary~\ref{tangentline},
the tangent line at an $\mathbb{F}_q$-point is a $1$-line or $2$-line, and vice versa. 
Since $\deg C=3$, there are no bi-tangents.
Hence
\begin{equation}\label{mutan}
\mu_1 + \mu_2 =9.
\end{equation}
Since $|p^{-1}(l)| =\delta$ if $l$ is a $\delta$-line,
\begin{equation}\label{muweightsum}
\mu_1 + 2 \mu_2 + 3 \mu_3= |\mathscr{I}| = 45.
\end{equation}
Additionally, since the total number of $\mathbb{F}_q$-lines is $21$,
\begin{equation}\label{musum}
\mu_0 + \mu_1 + \mu_2 + \mu_3 = 21.
\end{equation}

{\em Step}~2.
Suppose that $\mu_1=0$.
From (\ref{mutan}), (\ref{muweightsum}), (\ref{musum}),
we have
$
\mu_0 =3, \mu_2 = \mu_3 =9,
$
which is the case (ii).

{\em Step}~3.
Suppose that $\mu_1 \neq 0$.
Since (\ref{mutan}) and (\ref{muweightsum}),
$\mu_1 \equiv 0 \bmod 3$.
Hence there are at least three $1$-lines,
and hence there are at least three inflection $\mathbb{F}_4$-points.
Choose two inflection $\mathbb{F}_4$-points $Q_1$ and $Q_2$, and consider the line $l_0$ passing through these two points,
which is an $\mathbb{F}_4$-line.
Hence $l_0$ meets $C$ at another point $Q_0$,
which is also an $\mathbb{F}_4$-point.

{\em Claim}~1.
$Q_0$ is also a flex.

We need more notation. The linear equivalence relation of divisors on $C$ will be denoted by $\sim$, and a general line section on $C$ by $L$.
Here a general line section means a representative of the divisor cut out by a line on $C$, which makes sense up to the relation $\sim$.

\noindent \mbox{\em Proof of claim}~1.
Since
$Q_0 + Q_1 + Q_2 \sim L$ and $3Q_i \sim L$ for $i=1$ and $2$,
we have $3Q_0 \sim 3L -3Q_1 -3Q_2 \sim L$, which means that
$Q_0$ is a flex.
\qed

Hence the following property holds.
\begin{itemize}
\item[$(\dag)$] There are exactly three $\mathbb{F}_4$-lines passing through $Q_0$ besides $l_0$ and $T_{Q_0}(C)$, say $l_1, l_2, l_3$.
Each $l_i$ is a $3$-line.
\end{itemize}
Actually, since
\[
C(\mathbb{F}_4) = \{Q_0, Q_1, Q_2\} \sqcup \left(\sqcup_{i=1}^3 (l_i \cap C(\mathbb{F}_4) \setminus \{Q_0\}  )\right)
\]
and $|l_i \cap C(\mathbb{F}_4) \setminus \{Q_0\}| \leq 2$,
each $l_i$ is a $3$-line.

The six points of $C(\mathbb{F}_4)\setminus \{Q_0, Q_1, Q_2\}$
are named $\{P_i^{(j)} \mid i=1,2,3; j=1,2\}$
so that
$l_i \cap C(\mathbb{F}_4) =\{ Q_0, P_i^{(1)}, P_i^{(2)}\}.$

{\em Claim}~2.
$
\sum_{i=1}^3 ( P_i^{(1)}+ P_i^{(2)}) \sim 2L.
$

\noindent \mbox{\em Proof of claim}~2.
Since $Q_0 +  P_i^{(1)}+ P_i^{(2)} \sim L$ and $3Q_0 \sim L$,
we get $L + \sum_{i=1}^3 ( P_i^{(1)}+ P_i^{(2)}) \sim 3L$.
\qed

Since a nonsingular plane curve is projectively normal,
the divisor $\sum_{i=1}^3 ( P_i^{(1)}+ P_i^{(2)})$ on $C$ is cut out by a quadratic curve.
Let $D$ be the quadratic curve passing through those six points.
Suppose that $D$ is absolutely irreducible.
Then $D$ has exactly five $\mathbb{F}_4$-points if it is defined over $\mathbb{F}_4$, or at most four $\mathbb{F}_4$-points if it is not defined over $\mathbb{F}_4$ because an $\mathbb{F}_4$-point of $D$ is a point of $D \cap F_4(D)$;
both are absurd.
Therefore $D$ is a union of two lines $m, m'$.
If a line is not defined over $\mathbb{F}_4$, then $F_4(m) =m'$ and $D$ has only one $\mathbb{F}_4$-point: also absured.
Hence this split occurs over $\mathbb{F}_4$.
Since $\deg C = 3$, those six points split into two groups; three of them lie on $m$ and the remaining three lie on $m'$, and $P_i^{(1)}$ and $P_i^{(2)}$ do not belong the same group.
Hence we may assume that
$P_1^{(1)}, P_2^{(1)}, P_3^{(1)} \in m$ and
$P_1^{(2)}, P_2^{(2)}, P_3^{(2)} \in m'$.
Note that $m$ and $m'$ do not contain $Q_0$ nor $Q_1$ nor $Q_2$.

Apply the same arguments to $Q_1$ instead of $Q_0$ after $(\dag)$.
Since $Q_1$ does not lie on $m$ nor $m'$, there is a permutation
$(\sigma(1), \sigma(2), \sigma(3))$ of $(1,2,3)$ such that
$Q_1, P_i^{(1)}, P_{\sigma(i)}^{(2)}$ are collinear for $i=1, 2, 3$.
Similarly, there is another permutation $\tau$
such that
$Q_2, P_i^{(1)}, P_{\tau(i)}^{(2)}$ are collinear for $i=1, 2, 3$.
Therefore
\begin{equation}\label{LpassP11}
\left.
\begin{array}{ccc}
Q_0 +P_1^{(1)} + P_1^{(2)} & \sim & L\\
Q_1 +P_1^{(1)} + P_{\sigma(1)}^{(2)} & \sim & L\\
Q_2 +P_1^{(1)} + P_{\tau(1)}^{(2)} & \sim & L
\end{array}
\right\}
\end{equation}

{\em Claim}~3.
$\{ \sigma(1), \tau(1) \} = \{2, 3\}$.

\noindent \mbox{\em Proof of claim}~3.
If not, two of $\{ P_1^{(2)}, P_{\sigma(1)}^{(2)},  P_{\tau(1)}^{(2)}\}$
coincide.
For example, if $ P_1^{(2)}=P_{\sigma(1)}^{(2)}$,
then $Q_0, P_1^{(1)}, P_1^{(2)}=P_{\sigma(1)}^{(2)}, Q_1$
are collinear, which is impossible
because the line joining $Q_0$ and $Q_1$ is $l_0$.
Other cases can be handled by similar way.
\qed

By this claim,
\begin{equation}\label{twopermutations}
P_1^{(2)}+ P_{\sigma(1)}^{(2)} +P_{\tau(1)}^{(2)}
\sim L.
\end{equation}
Hence adding all equivalence relations in (\ref{LpassP11}),
together with (\ref{twopermutations}) we have
$
3 P_1^{(1)} +2L \sim 3L,
$
which implies $3P_1^{(1)} \sim L$.
Hence $ P_1^{(1)}$ is a flex.
Similarly we have that any $P_i^{(j)}$ is a flex.
Hence $\mu_1=9$.
Hence, from (\ref{mutan}), (\ref{muweightsum}) and (\ref{musum}) in Step~1,
$\mu_0=0$, $\mu_2=0$ and $\mu_3 =12$.
\qed

\begin{remark}
In Step~3 of the proof of Lemma~\ref{essentiallemma},
what we have shown is essentially that if a point of $C(\mathbb{F}_4)$ is flex,
then so are all points of  $C(\mathbb{F}_4)$.
If $C(\mathbb{F}_4)$ contains a flex, then $C$ is defined over $\mathbb{F}_4$ as an elliptic curve.
A sophisticated proof for the above fact may be possible by using the Jacobian variety, which coincides with the elliptic curve $C$.
For details, see the first part of \cite{ruc-sti1994}.
\end{remark}

\noindent \mbox{\em Proof of Theorem}~\ref{primarytheorem}.
When the case (ii) in Lemma~\ref{essentiallemma} occurs,
three $0$-lines are not concurrent;
Actually if three $0$-lines are concurrent, there is an $\mathbb{F}_4$-point
$Q$ outside $C$, which these $\mathbb{F}_4$-lines pass through.
The remaining two $\mathbb{F}_4$-lines pass through $Q$ can't cover all the points of $C(\mathbb{F}_4)$.

Hence we may choose coordinates $x_0, x_1, x_2$ so that
those $0$-lines are $\{x_0=0\}$, $\{x_1=0\}$ and $\{x_2=0\}$.
Since
$|\mathbb{P}^2(\mathbb{F}_4) \setminus \cup_{i=0}^2 \{x_i=0\}| = 9 = |C(\mathbb{F}_4)|$,
$C \in \mathscr{S}_4$
by Proposition~\ref{prop_characterization}.
Furthermore since $|\mathscr{S}_4|=1$ by Theorem~\ref{classification} (III-ii),
and $C_{(1,\omega, \omega^2)} \in \mathscr{S}_4$,
$C$ is projectively equivalent to 
to the curve
\[
x_0^3 + \omega x_1^3 + \omega^2 x_2^3 = 0.
\]

Next we consider the case (i) in Lemma~\ref{essentiallemma}.
In this case $C$ has the following properties:
\begin{enumerate}[(1)]
\item $C$ is nonsingular of degree $3$ defined over $\mathbb{F}_4$ with nine $\mathbb{F}_4$-points;
\item for any $P \in C(\mathbb{F}_4)$,
$i(T_P(C).C; P) =3$;
\item each point of $\mathbb{P}^2(\mathbb{F}_4) \setminus C(\mathbb{F}_4)$ lies on three tangent lines.
\end{enumerate}
Here we will confirm the property (3). Among the five $\mathbb{F}_4$-lines passing through $Q \in \mathbb{P}^2(\mathbb{F}_4) \setminus C(\mathbb{F}_4)$,
$\mu_{\delta}(Q)$ denotes the number of $\delta$-lines.
Since $\delta$ is either $1$ or $3$,
$\mu_1(Q) + 3 \mu_3(Q) = 9$ and $\mu_1(Q) + \mu_3(Q) =5.$
Hence $\mu_1(Q) =3$.

The proof of \cite[Lemma~7]{hir-sto-tha-vol1991} works well
under those three assumptions (1), (2), (3) for $C$.
To adapt their proof to our case, beware of a difference of notation;
their $q$ is our $\sqrt{q}$.
\qed

\section{Comparison of two maximal curves of degree $3$ over $\mathbb{F}_4$}
Lastly we compare two maximal curves of degree $3$
\[
C : \, x_0^3 + x_1^3 + x_2^3 = 0
\]
and
\[
D : \, x_0^3 + \omega x_1^3 + \omega^2 x_2^3 =0
\]
over $\mathbb{F}_4=\mathbb{F}_2[\omega]$.

Apparently, $C$ and $D$ are projectively equivalent over $\mathbb{F}_{2^6}$,
but not over $\mathbb{F}_{2^2}$
as we have seen.
We will show the function fields $\mathbb{F}_4(C)$
and $\mathbb{F}_4(D)$ are isomorphic over $\mathbb{F}_4$.
This is already guaranteed theoretically by R\"{u}ck and Stichtenoth \cite{ruc-sti1994}.
Here we will give an explicit isomorphism between those two fields.

Let $x = \frac{x_0}{x_2}|C$ and $y = \frac{x_1}{x_2}|C$.
Obviously $\mathbb{F}_4(C) =\mathbb{F}_4(x,y)$
with $x^3 + y^3 +1=0$.

\begin{theorem}
Three functions
\begin{align}
u &= 1 + \frac{x}{y+1} + \frac{1}{x+y+1} \notag\\
v &= \omega^2 \frac{x}{y+1} + \frac{1}{x+y+1} \label{uvw}\\
w &= \omega \frac{x}{y+1} + \frac{1}{x+y+1}\notag
\end{align}
satisfy
\[
u^3 + \omega v^3 + \omega^2 w^3=0.
\]
\end{theorem}
\proof
By straightforward  computation,
we have
\begin{align*}
(&(y+1)(x+y+1)w)^3\\
 =& ( \omega x(x+y+1) + (y+1) )^3\\
 =& 
 x^3(x+y+1)^3 + \omega^2 x^2(x+y+1)^2(y+1)
   + \omega x(x+y+1)(y+1)^2 + (y+1)^3,
 \end{align*}
\begin{align*}
(&(y+1)(x+y+1)v)^3\\
 =& ( \omega^2 x(x+y+1) + (y+1) )^3\\
 =& 
 x^3(x+y+1)^3 + \omega x^2(x+y+1)^2(y+1)
   + \omega^2 x(x+y+1)(y+1)^2 + (y+1)^3,
 \end{align*}
and
\begin{align*}
(&
(y+1)(x+y+1)u
)^3\\
 =& ( (y+1)(x+y+1)+ x(x+y+1) + (y+1) )^3  = g+h,\\
\end{align*}
where
\begin{align*}
g = & (y+1)^3(x+y+1)^3 + (y+1)^2(x+y+1)^2(x(x+y+1)+(y+1))\\
 &+ (y+1)(x+y+1)(x(x+y+1)+(y+1))^2,
\end{align*}
\begin{align*}
h= & (x(x+y+1)+(y+1))^3 \\
 =& x^3(x+y+1)^3 +x^2(x+y+1)^2(y+1) +
      x(x+y+1)(y+1)^2 +(y+1)^3.
\end{align*}
Hence
\begin{align*}
\omega^2((y+1)&(x+y+1)w)^3 + \omega((y+1)(x+y+1)v)^3 + h\\
=& (\omega^2 + \omega +1)x^3(x+y+1)^3 \\
 &+ (\omega^4 + \omega^2 +1)x^2(x+y+1)^2(y+1)\\
 &+ (\omega^3 + \omega^3 +1)x(x+y+1)(y+1)^2\\
 &+ (\omega^2 + \omega +1)(y+1)^3\\
 =& x(x+y+1)(y+1)^2.
\end{align*}
Therefore
\begin{align}
\omega^2&((y+1)(x+y+1)w)^3 + \omega((y+1)(x+y+1)v)^3 + ((y+1)(x+y+1)u)^3  
     \label{threeterms}\\
=& g + x(x+y+1)(y+1)^2 \notag \\
=&(y+1)(x+y+1)
\Bigl\{
(y+1)^2(x+y+1)^2 + x(y+1)(x+y+1)^2 \notag \\
&+ (y+1)^2(x+y+1) + x^2(x+y+1)^2
   + (y+1)^2 + x(y+1) \notag
\Bigl\}.
\end{align}
Since the sum of last two terms in the braces is
$(x+y+1)(y+1)$, $(x+y+1)$ divides the polynomial in the braces.
Hence (\ref{threeterms}) is equal to
\[
(y+1)^3(x+y+1)^3(\omega^2w^3 + \omega v^3 + u^3)
= (y+1)(x+y+1)^2f,
\]
where
\[
f =  (y+1)^2(x+y+1)+ x(y+1)(x+y+1)
+ (y+1)^2 + x^2(x+y+1)
   + (y+1)
\]
Continue the computation a little more:
\begin{align*}
f &= x(y+1)^2 + (y+1)^3 + x^2(y+1)+x(y+1)^2 
  + (y+1)^2 + x^3 + x^2(y+1) + (y+1)\\
 &= (y+1)^3 +(y+1)^2 + (y+1) + x^3 \\
  &= y^3 +x^3 +1 =0.
\end{align*}
As a conclusion, we have
$
u^3 + \omega v^3 + \omega^2 w^3 = 0.
$
\qed

\begin{corollary}
$\mathbb{F}_4(C) \cong \mathbb{F}_4(D).$
\end{corollary}
\proof
Trivially $\mathbb{F}_4(C) =\mathbb{F}_4(x,y)
 = \mathbb{F}_4(\frac{x}{y+1}, \frac{1}{x+y+1})$.
On the other hand, by definition of $u, v, w$ (\ref{uvw})
\[
\omega^2 \frac{v}{u} + \omega \frac{w}{u} =1- \frac{1}{u}.
\]
Hence $\mathbb{F}_4(D) \cong \mathbb{F}_4(\frac{v}{u}, \frac{w}{u})
 = \mathbb{F}_4(u,v,w)$.
Since
\[
\begin{pmatrix}
u \\ v\\ w
\end{pmatrix}
=
\begin{pmatrix}
1 & 1& 1\\
0 & \omega^2 & 1 \\
0 & \omega & 1
\end{pmatrix}
\begin{pmatrix}
1 \\ \frac{x}{y+1} \\\frac{1}{x+y+1}
\end{pmatrix},
\]
we know 
$ \mathbb{F}_4(u,v,w)=\mathbb{F}_4(\frac{x}{y+1}, \frac{1}{x+y+1}).$
Summing up, we get
$\mathbb{F}_4(D) \cong \mathbb{F}_4(C)$.
\qed

\end{document}